\documentclass[10pt,leqno,twoside]{amsart}

\usepackage{enumerate}

\usepackage{amssymb}

\usepackage{amssymb,amsthm,amsmath,eucal,mathrsfs,verbatim}

\usepackage{footnote}

\usepackage{tikz, subfigure, xcolor} 

\usepackage{pgfplots}

\usepackage{cite}

\usepackage{amsmath}

\usepackage{amsthm}

\usepackage{amssymb}

\usepackage{amsfonts}

\usepackage{dsfont}

\usepackage{hyperref}

\newtheorem{theorem}{Theorem}[section]

\newtheorem{lemma}[theorem]{Lemma}

\newtheorem{proposition}[theorem]{Proposition}

\theoremstyle{definition}

\numberwithin{equation}{section}

\newcommand{\dist}{\mathrm{dist}}      

\renewcommand{\Im}{{\ensuremath{\mathrm{Im\,}}}} 

\renewcommand{\Re}{{\ensuremath{\mathrm{Re\,}}}} 

\newcommand\restr[2]{{
  \left.\kern-\nulldelimiterspace 
  #1 
  \vphantom{\big|} 
  \right|_{#2} 
  }}

\usepackage{eucal}

\title[Increasing stability]{Increasing stability for the inverse problem for the Schr{\"o}dinger equation}
\author[Choudhury, Heck]{Anupam Pal Choudhury $^* $ , Horst Heck $^\dagger $ } 
\date{}

\keywords{Inverse problems, stability estimates, Schr{\"o}dinger equation \\
$^*$ Technische Universit{\"a}t Darmstadt, Institute of Mathematics, Schlossgartenstr. 7, 
64289 Darmstadt. Email: anupampcmath@gmail.com.  This author is supported by the DFG International Research Training Group IRTG 1529 on Mathematical Fluid Dynamics at TU Darmstadt \\ 
$^\dagger $ Departement Engineering and Information Technology, Bern University of Applied Sciences, Jlcoweg 1, CH-3400 Burgdorf, Switzerland.  
Email: horst.heck@bfh.ch  }

\begin{document}

\maketitle
\abstractname.{ In this article, we study the increasing stability property
for the determination of the potential in the Schr{\"o}dinger equation from
partial data. We shall assume that the inaccessible part of the boundary is
flat and homogeneous boundary condition is prescribed on this part. In
contrast to earlier works, we are able to deal with the case when potentials
have some Sobolev regularity and also need not be compactly supported inside the domain. }

\section{Introduction}
Let us consider the boundary value problem for the  Schr\"odinger equation 
\begin{equation}
(\Delta+k^{2}+q(x)) u(x)=0, \ \text{in} \ \Omega, 
\label{main-1}
\end{equation}
posed in a bounded domain $\Omega \subset \mathbb{R}^{3} $ with smooth boundary. The boundary data
\begin{equation}
u(x) = f(x) \ \text{on} \ \partial \Omega.
\label{main-2}
\end{equation}
is assumed to be of the class $H^{\frac{1}{2}}(\partial \Omega) $, and $q$ is
real-valued and satisfies $q \in H^{s}(\Omega)  $, for some $s > \frac{3}{2} $. 
Note, that, by Sobolev embedding, this yields, that the potentials are in fact
H{\"o}lder continuous.
Without loss of generality, we shall assume that the wave number $k \geq 1 $.\\
For $N > 0 $ and $s>\frac{3}{2} $, let us define the set of potentials 
\begin{equation}
\mathcal{Q}_{N}:= \{q: \Vert q \Vert_{H^{s}(\Omega)} \leq N \}.
\notag
\end{equation}
In this article, we shall consider a bounded domain $\Omega $ with smooth boundary such that $\Omega \subset \{x: x_{3} < 0 \} $ and a part of the boundary $\Gamma_{0} $ (which we shall also refer as the inaccessible part of the boundary) is contained in the plane $\{x: x_{3}=0 \} $. We shall assume that the support of $f$ is contained in $\Gamma:=\partial \Omega \setminus \Gamma_{0} $.\\
Let 
\begin{equation}
\mathcal{C}_{q}:= \Big\{\Big(u\Big\vert_{\Gamma},\frac{\partial u}{\partial \nu}\Big\vert_{\Gamma} \Big), \text{where $u$ is a solution to \eqref{main-1} and $u=0 $ on $\Gamma_{0} $} \Big\}
\notag
\end{equation}
denote the partial Cauchy data and $\frac{\partial u}{\partial \nu}\vert_{\Gamma} \in H^{-\frac{1}{2}}(\Gamma) $.
We can define a distance in the set of partial Cauchy data as 
\begin{equation}
\dist(\mathcal{C}_{q_{1}},\mathcal{C}_{q_{2}}):= \text{max} \Big\{\max_{(f,g) \in \mathcal{C}_{q_{1}}} \min_{(\tilde{f},\tilde{g}) \in \mathcal{C}_{q_{2}}} \frac{\Vert (f,g)-(\tilde{f},\tilde{g}) \Vert_{H^{\frac{1}{2}} \oplus H^{-\frac{1}{2}}}}{\Vert (f,g) \Vert_{H^{\frac{1}{2}} \oplus H^{-\frac{1}{2}}}}, \max_{(f,g) \in \mathcal{C}_{q_{2}} } \min_{(\tilde{f},\tilde{g}) \in \mathcal{C}_{q_{1}}} \frac{\Vert (f,g)-(\tilde{f},\tilde{g}) \Vert_{H^{\frac{1}{2}} \oplus H^{-\frac{1}{2}}}}{\Vert (f,g) \Vert_{H^{\frac{1}{2}} \oplus H^{-\frac{1}{2}}}} \Big\},
\notag
\end{equation}
where $\Vert (f,g) \Vert_{H^{\frac{1}{2}} \oplus H^{-\frac{1}{2}}} = (\Vert f \Vert^{2}_{H^{\frac{1}{2}}(\Gamma)}+\Vert g \Vert^{2}_{H^{-\frac{1}{2}}(\Gamma)})^{\frac{1}{2}} $.\\
Our aim, here, is to address the question of stability of the recovery of the potential $q$ from the knowledge of the partial Cauchy data $\mathcal{C}_{q}$ and to study the behaviour of the stability estimates as the wave number $k$ grows. The unique identification of the potential $q$ from $\mathcal{C}_{q} $ was established in the work \cite{Isakov-1}. \\
Starting with the work \cite{Calderon-InverseProblemPaper} and following the
impetus provided by the work \cite{Sylvester-Uhlmann-CalderonProblemPaper},
such problems started receiving intense consideration. The question of
stability in the case of full data (and $k=0$) was considered in
\cite{Alessandrini-StabilityPaper} and a logarithmic stability estimate was
established. It was also shown that this is the optimal result one can achieve. In the partial data case (with $k=0$), a double logarithmic type stability estimate was established in \cite{Heck-Wang-StabilityPaper} following the work \cite{Bukhgeim-Uhlmann-InverseProblemPaper} which dealt with the issue of unique identification. We would also like to refer to the work \cite{Heck-StabilityPaper} in this context. In the case of domains under consideration (with $k=0$), it was shown in \cite{Heck-Wang-Optimal} that a logarithmic type stability estimate can be established even from partial data. \\
In order to improve the logarithmic type stability estimates (which
means that the problem is severely ill-posed and therefore
inconvenient also from a numerical point of view) to Lipschitz-type stability
estimates, the corresponding problem with $k\not=0$ started receiving attention. It was found in many works (see \cite{Isakov-2,INUW, IW, Liang, NUW}) and in the context of different models that a growing $k$ tends to improve the stability, a property which was termed as increasing stability. In this article, we shall investigate this property in case of the domains stated above and endeavour to improve the logarithmic stability estimate established in \cite{Heck-Wang-Optimal}. \\  
We would like to remark that the property of increasing stability in similar
domains was also studied in \cite{Liang}. In that article, the author assumed
the condition $\frac{\partial u}{\partial \nu}=0 \ \text{on}\ \Gamma_{0} $
instead of $u=0 \ \text{on} \ \Gamma_{0} $ that we have assumed here.
Nevertheless our proof with minor modifications (see \cite{Isakov-1}) would
also hold true in that case. Moreover, here we assume only Sobolev
regularity of the potentials in contrast to the assumption of potentials in $C^{1}(\Omega) $ considered in \cite{Liang}. We also do not assume that the difference of the potentials vanishes near the boundary $\partial \Omega $.\\
Our main result on the stability of recovery of the potential $q$ from the Cauchy data $\mathcal{C}_{q} $ reads as follows.
\begin{theorem}\label{main-theorem}
Let $\Omega \subset \mathbb{R}^{3} $ be a bounded domain as described above. Also assume that $R>0 $ be a large real number such that $\Omega \subset B(0,R) $. Let $\mathcal{C}_{q_{1}}, \mathcal{C}_{q_{2}} $ denote the partial Cauchy data corresponding to the potentials $q_{1}, q_{2} \in \mathcal{Q}_{N} $ respectively. Then there exist constants $C, \tilde{\alpha}, \eta > 0 $ such that
\begin{equation}
\Vert q_{1}-q_{2} \Vert_{L^{\infty}(\Omega)} \leq C \Big(e^{6Rk} \ \dist(\mathcal{C}_{q_{1}},\mathcal{C}_{q_{2}})+ \frac{1}{(k+\frac{E}{5R})^{\tilde{\alpha}}} \Big)^{\frac{\eta}{2(1+s)}},
\label{main-10}
\end{equation}
where $E= \vert \log \ \dist(\mathcal{C}_{q_{1}},\mathcal{C}_{q_{2}}) \vert $. The constant $C$ depends on $\Omega, N \ \text{and}\ s $ only and the constants $\tilde{\alpha},\eta $ depend on $s$ only.
\end{theorem}
It can be observed from \eqref{main-10} that as the positive constant $k$ grows, the term in the right-hand side with the logarithmic part in the denominator decays to zero and the first term (the Lipschitz part) dominates. Thus the logarithmic stability is improved to a Lipschitz-type stability estimate exhibiting the property of increasing stability.\\  
The above result should also hold true, with minor modifications, for any dimension $n>3$. To simplify the presentation in terms of the CGO solutions, we have restricted ourselves to the case $n=3$.  \\
In Section 2, we recollect some preliminary results that shall be necessary in the proof of the stability estimates. In Section 3, we introduce appropriate solutions to \eqref{main-1} and proceed to derive the desired stability estimates. 
\section{Some preliminary results}
In this section, we recollect some preliminary results which we shall use later in the proofs. We begin by stating a result on the existence of complex geometric optics (CGO) solutions to \eqref{main-1}.
\begin{lemma}(see \cite{INUW,Sylvester-Uhlmann-CalderonProblemPaper})\label{CGO-existence}
Let $s>\frac{3}{2} $. Assume that $\zeta= \Re  \zeta+i \ \Im  \zeta $
satisfies $\vert \Re \zeta \vert^{2}=k^{2}+\vert \Im \zeta \vert^{2} $ and
$\Re \zeta \cdot \Im \zeta=0 $, that is, $\zeta \cdot \zeta=k^{2} $. Then there 
exist constants $C_{*}$ and $C>0 $, independent of $k$, such that if $\vert \Im \zeta \vert > C_{*} \Vert q \Vert_{H^{s}(\Omega)}$, then there exists a solution $u$ to \eqref{main-1} of the form
\[u(x)=e^{i \zeta \cdot x}(1+\psi(x)) \] where \[\Vert \psi
\Vert_{H^{s}(\Omega)} \leq \frac{C}{\vert \Im\zeta \vert} \Vert q \Vert_{H^{s}(\Omega)} .\] 
\end{lemma}
In the next section, we shall choose $\zeta $ suitably so as to be able to use the above lemma to infer the existence of CGO solutions with the error terms satisfying the above estimates. We shall also need the following Green's identity which can be proved following 
\cite{Alessandrini-StabilityPaper, INUW}. 
\begin{proposition}
Let $u_{j} $ and $\mathcal{C}_{q_{j}} $ be solution and Cauchy data for the equation \eqref{main-1} corresponding to the potential $q_{j} \ (j=1,2) $. Then 
\begin{equation}
\Big\vert \int_{\Omega} (q_{2}-q_{1})\ u_{1} u_{2}\ dx \Big\vert \leq \Big\Vert \Big(u_{1}, \frac{\partial u_{1}}{\partial \nu} \Big) \Big\Vert_{H^{\frac{1}{2}}\oplus H^{-\frac{1}{2}}} 
\Big\Vert \Big(u_{2}, \frac{\partial u_{2}}{\partial \nu} \Big)
\Big\Vert_{H^{\frac{1}{2}}\oplus H^{-\frac{1}{2}}} \ \dist(\mathcal{C}_{q_{1}}, \mathcal{C}_{q_{2}}).
\notag
\end{equation}
\end{proposition}
Using the equation \eqref{main-1}, it can be proved (see \cite{FSU,INUW}) that 
\begin{equation}
\begin{aligned}
\Big\Vert \Big(u_{l}, \frac{\partial u_{l}}{\partial \nu} \Big) \Big\Vert_{H^{\frac{1}{2}}\oplus H^{-\frac{1}{2}}} &\leq C k^{2} \Vert u_{l} \Vert_{L^{2}(\Omega)} + C \Vert \nabla u_{l} \Vert_{L^{2}(\Omega)} \\
& \leq C k^{2} \Vert u_{l} \Vert_{H^{1}(\Omega)}.
\end{aligned}
\notag
\end{equation}
Using this together with the above proposition, we can derive 
\begin{equation}
\Big\vert \int_{\Omega} (q_{2}-q_{1})\ u_{1} u_{2}\ dx \Big\vert \leq C k^{4}
\Vert u_{1} \Vert_{H^{1}(\Omega)} \Vert u_{2} \Vert_{H^{1}(\Omega)} \ \dist(\mathcal{C}_{q_{1}}, \mathcal{C}_{q_{2}}).
\label{prel-17}
\end{equation}
In what follows, we shall also require the following quantitative version of the Riemann-Lebesgue lemma. For the proofs of the results, we refer to \cite{Choudhury-Heck, Heck-Wang-Optimal}. 
\begin{lemma}
Let $\Omega \subset \mathbb{R}^{n} $ be a bounded domain with $C^{1}$ boundary and let $f \in C^{0,\alpha}(\overline{\Omega}) $ for some $\alpha \in (0,1)$. Let $\tilde{f}$ denote the extension of $f $ to $\mathbb{R}^{n}$ by zero. Then there exist $\tilde{\delta} > 0$ and $C>0$ such that
\begin{equation}
\Vert \tilde{f}(\cdot-y)-\tilde{f}(\cdot) \Vert_{L^{1}(R)^{n}} \leq C \vert y \vert^{\alpha},
\notag
\end{equation}
for any $y \in \mathbb{R}^{n} $ with $\vert y \vert < \tilde{\delta} $.
\end{lemma}
\begin{lemma}\label{estimate}
Let $f \in L^{1}(\mathbb{R}^{n})$ and suppose there exist constants $\tilde{\delta} > 0, \ C_{0}>0 $ and $\alpha \in (0,1) $ such that for $ \vert y \vert < \tilde{\delta}$,
\begin{equation}
\Vert f(\cdot-y)-f(y) \Vert_{L^{1}(\mathbb{R}^{n})} \leq C_{0} \vert y \vert^{\alpha}.
\label{prel18}
\end{equation}
Then there exist constants $C>0 $ and $\epsilon_{0}>0 $ such that for any $0<\epsilon<\epsilon_{0} $, we have the inequality
\begin{equation}
\vert \mathcal{F} f (\xi) \vert \leq C (e^{-\frac{\epsilon^{2} \vert \xi \vert^{2}}{4 \pi}}+ \epsilon^{\alpha} ),
\label{prel19}
\end{equation}
where the constant $C$ depends on $C_{0}, \Vert f \Vert_{L^{1}}, n , \tilde{\delta} $ and $\alpha $.
\end{lemma}
By assumption, the potentials $q \in H^{s}(\Omega)$ with $s>\frac{3}{2}$ and therefore there exists $\alpha>0$ such that $q \in C^{0,\alpha}(\overline{\Omega})$. The conclusions of the Lemma \ref{estimate}, therefore,  hold true for the potentials $q$. 
\section{CGO and the stability estimates}
In this section, we shall construct appropriate solutions to \eqref{main-1} via CGO solutions as described in Lemma \ref{CGO-existence}. In order to do so, we introduce a change of coordinates as follows (see also \cite{Choudhury-Heck, Heck-Wang-Optimal, Isakov-1, Liang}).\\
Given $\xi=(\xi_{1},\xi_{2},\xi_{3}) \in \mathbb{R}^{3} $, the new coordinate representation 
is obtained by rotating the standard axes in a manner such that under the transformed coordinates, the representation of $\xi$, which we shall denote henceforth by $\tilde{\xi} $, is of the form 
\begin{equation}
\tilde{\xi}=(\tilde{\xi}_{1},0,\tilde{\xi}_{3}), \ \text{where}\ \tilde{\xi}_{1}=(\xi_{1}^{2}+\xi_{2}^{2})^{\frac{1}{2}} \ \text{and} \ \tilde{\xi}_{3}=\xi_{3}.
\notag
\end{equation}
Let $\tilde{x} $ denote the representation of $x $ in this new coordinates. It is easy to see that for $x,y \in \mathbb{R}^{3} $, $\sum_{i=1}^{3} x_{i} \cdot y_{i}= \sum_{i=1}^{3} \tilde{x}_{i} \cdot \tilde{y}_{i} $.\\
In the transformed coordinates, let us choose
\begin{equation}
\begin{aligned}
\tilde{\zeta}_{1}&= \Big(-\frac{\tilde{\xi}_{1}}{2}+ \tau \tilde{\xi}_{3} , -i \Big(\vert \xi \vert^{2}(\frac{1}{4}+\tau^{2})-k^{2} \Big)^{\frac{1}{2}}, -\frac{\tilde{\xi}_{3}}{2}-\tau \tilde{\xi}_{1}  \Big), \\ 
\tilde{\zeta}_{2}&= \Big( -\frac{\tilde{\xi}_{1}}{2}- \tau \tilde{\xi}_{3} , i \Big(\vert \xi \vert^{2}(\frac{1}{4}+\tau^{2})-k^{2} \Big)^{\frac{1}{2}}, -\frac{\tilde{\xi}_{3}}{2}+\tau \tilde{\xi}_{1} \Big), 
\end{aligned}
\label{estimate-1}
\end{equation}
where $\tau $ is a positive real number.
We note that in the original coordinates $\zeta_{1},\zeta_{2} $ are of the form
\begin{equation}
\begin{aligned}
&\zeta_{1}: \begin{cases}
\zeta_{1,1}=(-\frac{\tilde{\xi_{1}}}{2}+\tau \tilde{\xi}_{3}) \frac{\xi_{1}}{\sqrt{\xi_{1}^{2}+\xi_{2}^{2}}} + i \frac{\xi_{2}}{\sqrt{\xi_{1}^{2}+\xi_{2}^{2}}} (\vert \xi \vert^{2} (\frac{1}{4}+\tau^{2})-k^{2})^{\frac{1}{2}} \\
\zeta_{1,2}= (-\frac{\tilde{\xi_{1}}}{2}+\tau \tilde{\xi}_{3}) \frac{\xi_{2}}{\sqrt{\xi_{1}^{2}+\xi_{2}^{2}}} - i \frac{\xi_{1}}{\sqrt{\xi_{1}^{2}+\xi_{2}^{2}}} (\vert \xi \vert^{2} (\frac{1}{4}+\tau^{2})-k^{2})^{\frac{1}{2}}\\
\zeta_{1,3}=-\frac{\tilde{\xi}_{3}}{2}-\tau \tilde{\xi}_{1} = -\frac{\xi_{3}}{2} -\tau (\xi_{1}^{2}+\xi_{2}^{2})^{\frac{1}{2}} 
\end{cases} \\
&\zeta_{2}: \begin{cases}
\zeta_{2,1}= (-\frac{\tilde{\xi_{1}}}{2}-\tau \tilde{\xi}_{3}) \frac{\xi_{1}}{\sqrt{\xi_{1}^{2}+\xi_{2}^{2}}} - i \frac{\xi_{2}}{\sqrt{\xi_{1}^{2}+\xi_{2}^{2}}} (\vert \xi \vert^{2} (\frac{1}{4}+\tau^{2})-k^{2})^{\frac{1}{2}}\\
\zeta_{2,2}=(-\frac{\tilde{\xi_{1}}}{2}-\tau \tilde{\xi}_{3}) \frac{\xi_{2}}{\sqrt{\xi_{1}^{2}+\xi_{2}^{2}}} + i \frac{\xi_{1}}{\sqrt{\xi_{1}^{2}+\xi_{2}^{2}}} (\vert \xi \vert^{2} (\frac{1}{4}+\tau^{2})-k^{2})^{\frac{1}{2}} \\
\zeta_{2,3}= -\frac{\tilde{\xi}_{3}}{2}+\tau \tilde{\xi}_{1} = -\frac{\xi_{3}}{2} +\tau (\xi_{1}^{2}+\xi_{2}^{2})^{\frac{1}{2}} 
\end{cases}
\end{aligned}
\label{estimate-2}
\end{equation}
where $\zeta_{i,j} $ denote the j-th coordinate of $\zeta_{i}$. 
Let us also define the reflections of $\tilde\zeta_i$ on the plane $\xi_3=0$,
i.e.
\begin{equation}
\begin{aligned}
\tilde{\zeta}^{*}_{1}&= \Big(-\frac{\tilde{\xi}_{1}}{2}+ \tau \tilde{\xi}_{3} , -i \Big(\vert \xi \vert^{2}(\frac{1}{4}+\tau^{2})-k^{2} \Big)^{\frac{1}{2}}, \frac{\tilde{\xi}_{3}}{2}+\tau \tilde{\xi}_{1}  \Big),\\
\tilde{\zeta}^{*}_{2}&= \Big( -\frac{\tilde{\xi}_{1}}{2}- \tau \tilde{\xi}_{3} , i \Big(\vert \xi \vert^{2}(\frac{1}{4}+\tau^{2})-k^{2} \Big)^{\frac{1}{2}}, \frac{\tilde{\xi}_{3}}{2}-\tau \tilde{\xi}_{1} \Big) .
\end{aligned}
\label{estimate-3}
\end{equation}
 It is easy to see from \eqref{estimate-1}--\eqref{estimate-3} that for $j=1,2$,
\begin{equation}
\vert \Re\ \zeta_{j} \vert^{2}=\vert \xi \vert^{2} (\frac{1}{4}+\tau^{2}) , \
\text{and}\ 
\vert \Im\ \zeta_{j} \vert^{2}=\vert \xi \vert^{2} (\frac{1}{4}+\tau^{2})-k^{2},
\notag
\end{equation} 
and $\zeta_{1}+\zeta_{2}, \zeta^{*}_{1}+\zeta^{*}_{2},\zeta_{1}+\zeta^{*}_{2},\zeta^{*}_{1	}+\zeta_{2} $ are real vectors. Also $\zeta_{1} \cdot \zeta_{1}=\zeta_{2} \cdot \zeta_{2}=\zeta_{1}^{*} \cdot \zeta_{1}^{*}=\zeta_{2}^{*} \cdot \zeta_{2}^{*}  =k^{2}$.\\
We next extend the potentials $q_{1}, q_{2} $ onto the whole of
$\mathbb{R}^{3} $ as even functions in $x_{3} $. Lemma \ref{CGO-existence}
then guarantees the existence of CGO solutions to the equation \eqref{main-1}
in $\mathbb{R}^{3} $ of the form $e^{i \zeta_{j} \cdot x}(1+ w_{j}) $ and
$e^{i \zeta^{*}_{j} \cdot x}(1+w^{*}_{j}) $ for $j=1,2 $ with the remainder
terms satisfying the estimate $\Vert w_{j} \Vert_{H^{s}(\Omega)} \leq
\frac{C}{\vert \Im  \ \zeta_{j} \vert} \Vert q_{j} \Vert_{H^{s}(\Omega)} $.\\
Let us define
\begin{equation}
\begin{aligned}
u_{1}(x)&= e^{i \zeta_{1} \cdot x}(1+ w_{1}) - e^{i \zeta^{*}_{1} \cdot x}(1+w^{*}_{1}), \
u_{2}(x)= e^{i \zeta_{2} \cdot x}(1+ w_{2}) - e^{i \zeta^{*}_{2} \cdot x}(1+w^{*}_{2}).
\end{aligned}
\label{estimate-4}
\end{equation}
It can be easily checked from the definitions \eqref{estimate-4} that the functions $u_{j}  \ (j=1,2) $ satisfy the equations \eqref{main-1} in $\mathbb{R}^{3}_{-} $ with potentials $q_{1},q_{2} $ respectively and $u_{j}(x)=0 $ on $x_{3}=0 $.\\
With all this preparation in place, we now proceed to derive the stability estimates. 
\subsection{Derivation of the stability estimates}
Let us denote $M= C_{*} N $. Then provided $\vert \Im \ \zeta_{j} \vert > M $, the estimate
\begin{equation}
\Vert w_{j} \Vert_{H^{s}(\Omega)} \leq \frac{C}{\vert \Im \ \zeta_{j}\vert}
\Vert q_{j} \Vert_{H^{s}(\Omega)} \leq \frac{C N}{\vert \Im \ \zeta_{j} \vert} \leq \frac{C  N}{C_{*} N} = C,
\notag
\end{equation}
holds true.
Let $\Omega \subset B(0,R) $ for a fixed $R (>> 1)$ large enough. Then since $\vert e^{i \zeta_{j} \cdot x}\vert \leq e^{\vert Im \ \zeta_{j} \vert \vert x \vert} $, we can write 
\begin{equation}
\begin{aligned}
\Vert u_{j} \Vert_{H^{1}(\Omega)} &\leq 2 e^{R [\vert \xi \vert^{2}(\frac{1}{4}+\tau^{2})-k^{2} ]^{\frac{1}{2}}} \Vert 1+ w_{j} \Vert_{H^{s}(\Omega)} \leq
 C e^{R [\vert \xi \vert^{2} (\frac{1}{4}+\tau^{2})-k^{2} ]^{\frac{1}{2}}},
\end{aligned}
\notag
\end{equation}
since $s > \frac{3}{2} >1 $. 
Using this in \eqref{prel-17}, we see that
\begin{equation}
\Big\vert \int_{\Omega} (q_{2}-q_{1}) u_{1} u_{2} \ dx \Big\vert \leq C k^{4}
e^{2 R \Big[\vert \xi \vert^{2} (\frac{1}{4}+\tau^{2})-k^{2}
\Big]^{\frac{1}{2}}}\ \dist(\mathcal{C}_{q_{1}},\mathcal{C}_{q_{2}}),
\label{estimate-20}
\end{equation} 
provided $\vert \Im \ \zeta_{j} \vert > M $, that is, $ \vert \xi \vert^{2} \Big(\frac{1}{4}+\tau^{2} \Big) > M^{2}+k^{2} $.\\
Let us denote $q_{0}=q_{2}-q_{1} $. Using the definitions of $u_{1}, u_{2} $ from \eqref{estimate-4}, we can write 
\begin{equation}
\begin{aligned}
\int_{\Omega} q_{0} u_{1} u_{2} \ dx &= \int_{\Omega} q_{0}(x) \Big[e^{i (\zeta_{1}+\zeta_{2}) \cdot x} (1+w_{1})(1+w_{2}) + e^{i (\zeta^{*}_{1}+\zeta^{*}_{2}) \cdot x}(1+ w_{1}^{*})(1+w_{2}^{*}) \\
& \qquad -e^{i (\zeta_{1}+\zeta^{*}_{2}) \cdot x}(1+w_{1})(1+w_{2}^{*}) -e^{i (\zeta^{*}_{1}+\zeta_{2}) \cdot x}(1+w^{*}_{1})(1+w_{2}) \Big] \ dx \\
&= \int_{\Omega} q_{0}(x) [e^{-i \xi \cdot x} +e^{-i \xi^{*} \cdot x} ] \ dx +\int_{\Omega} q_{0}(x) f(x, w_{1}, w_{2}, w_{1}^{*}, w_{2}^{*})\ dx \\
&\qquad -\int_{\Omega} q_{0}(x) [e^{i (\zeta_{1}+\zeta^{*}_{2}) \cdot x} + e^{i (\zeta^{*}_{1}+\zeta_{2}) \cdot x}] \ dx \\
&=\mathcal{F}q_{0}(\xi)+ \int_{\Omega} q_{0}(x) f(x, w_{1}, w_{2}, w_{1}^{*}, w_{2}^{*}) \ dx -\int_{\Omega} q_{0}(x) [e^{i (\zeta_{1}+\zeta^{*}_{2}) \cdot x} + e^{i (\zeta^{*}_{1}+\zeta_{2}) \cdot x}] \ dx,
\end{aligned}
\label{estimate-21}
\end{equation}
where 
\begin{equation}
\begin{aligned}
f = e^{-i \xi \cdot x} (w_{1}+w_{2}+w_{1}w_{2}) &+ e^{- i \xi^{*} \cdot x} (w_{1}^{*} + w_{2}^{*} + w_{1}^{*} w_{2}^{*}) -e^{i (\zeta^{*}_{1}+ \zeta_{2}) \cdot x} (w^{*}_{1}+w_{2}+w^{*}_{1} w_{2}) \\ &
-e^{i (\zeta_{1}+\zeta_{2}^{*}) \cdot x} (w_{1}+w_{2}^{*}+w_{1} w_{2}^{*}).
\end{aligned}
\label{estimate-22}
\end{equation}
Next since $\zeta_{1}+\zeta^{*}_{2}, \tilde{\zeta}_{1}+\tilde{\zeta}_{2}^{*}, \zeta^{*}_{1}+\zeta_{2} $ and $\tilde{\zeta}_{1}^{*}+\tilde{\zeta}_{2} $ are real vectors, we note that 
\begin{equation}
\begin{aligned}
&(\zeta_{1}+\zeta^{*}_{2}) \cdot x= (\tilde{\zeta}_{1}+\tilde{\zeta}_{2}^{*}) \cdot \tilde{x} =-\tilde{\xi}_{1} \tilde{x}_{1} - 2 \tau \tilde{\xi}_{1} \tilde{x}_{3} =- [\xi' \cdot x' + 2 \tau \vert \xi' \vert x_{3}], \\ 
&(\zeta^{*}_{1}+\zeta_{2}) \cdot x= (\tilde{\zeta}_{1}^{*}+\tilde{\zeta}_{2}) \cdot \tilde{x} =-\tilde{\xi}_{1} \tilde{x}_{1} + 2 \tau \tilde{\xi}_{1} \tilde{x}_{3} =- [\xi' \cdot x' - 2 \tau \vert \xi' \vert x_{3}],
\end{aligned}
\notag
\end{equation}
where $\xi'=(\xi_{1},\xi_{2}), x'=(x_{1},x_{2}) $, and therefore we can write
\begin{equation}
\begin{aligned}
&\int_{\Omega} q_{0}(x) e^{i (\zeta_{1}+\zeta^{*}_{2}) \cdot x} \ dx=
\mathcal{F}q_{0}(\xi',2 \tau \vert \xi' \vert) \ \text{and} \ \int_{\Omega} q_{0}(x) e^{i (\zeta^{*}_{1}+\zeta_{2}) \cdot x} \ dx= \mathcal{F}q_{0}(\xi',-2 \tau \vert \xi' \vert)  .
\end{aligned}
\label{estimate-23}
\end{equation}
Using the version of Riemann-Lebesgue Lemma stated in Lemma \ref{estimate}, the terms in \eqref{estimate-23} can be estimated as
\begin{equation}
\begin{aligned}
\vert \mathcal{F}q_{0}(\xi',2 \tau \vert \xi' \vert) \vert + \vert \mathcal{F}q_{0}(\xi',-2 \tau \vert \xi' \vert) \vert \leq C \Big[e^{-\frac{\epsilon^{2} (1+4 \tau^{2}) \vert \xi' \vert^{2}}{4 \pi}} + \epsilon^{\alpha} \Big], \ \text{where}\ \alpha \in (0,1),
\end{aligned}
\label{estimate-24}
\end{equation}
and for any $\epsilon< \epsilon_{0} $ with $\epsilon_{0} $ defined as in Lemma \ref{estimate}. \\
Also for $\vert \xi \vert^{2} (\frac{1}{4}+\tau^{2}) > M^{2}+k^{2} $, we can use the estimates for the remainder terms $w_{j} $ to derive
\begin{equation}
\vert \int_{\Omega} q_{0} f(x,w_{1},w_{2},w_{1}^{*}, w_{2}^{*}) \vert \ dx \leq C \Vert q_{0} \Vert_{L^{2}} \Vert f \Vert_{L^{2}} \leq \frac{C }{[\vert \xi \vert^{2} (\frac{1}{4}+ \tau^{2})- k^{2} ]^{\frac{1}{2}}}  .
\label{estimate-25}
\end{equation}
Combining the above estimates \eqref{estimate-20}--\eqref{estimate-25}, it follows that provided $\vert \xi \vert^{2} (\frac{1}{4}+\tau^{2}) > M^{2}+k^{2}$ holds, we have
 \begin{equation}
 \begin{aligned}
\vert \mathcal{F}q_{0}(\xi)\vert  &\leq C \Big[k^{4} e^{2R [\vert \xi
\vert^{2}(\frac{1}{4}+\tau^{2}) - k^{2}]^{\frac{1}{2}}} \dist(\mathcal{C}_{q_{1}},\mathcal{C}_{q_{2}}) + e^{-\frac{\epsilon^{2} (1+ 4 \tau^{2}) \vert \xi' \vert^{2}}{4 \pi}} + \epsilon^{\alpha} + \frac{1}{[\vert \xi \vert^{2}(\frac{1}{4}+\tau^{2})-k^{2}]^{\frac{1}{2}}}  \Big] \\
\Rightarrow \vert \mathcal{F}q_{0}(\xi) \vert^{2} &\leq C \Big[k^{8} e^{4R
[\vert \xi \vert^{2}(\frac{1}{4}+\tau^{2}) - k^{2}]^{\frac{1}{2}}} \dist(\mathcal{C}_{q_{1}},\mathcal{C}_{q_{2}})^{2} + e^{-\frac{\epsilon^{2} (1+ 4 \tau^{2}) \vert \xi' \vert^{2}}{2 \pi}} + \epsilon^{2 \alpha} + \frac{1}{\vert \xi \vert^{2}(\frac{1}{4}+\tau^{2})-k^{2}} \Big].
  \end{aligned}
 \label{estimate-26}
 \end{equation}
Our strategy next is to estimate the $H^{-1} $ norm of $q_{0}$ and then use
the interpolation inequality to derive an estimate for the $L^{\infty} $ norm of $q_{0}$. It will be worthwhile to note at this point that it is sufficient to derive the stability estimates when $dist(\mathcal{C}_{q_{1}},\mathcal{C}_{q_{2}}) < \delta $ for some $\delta \in (0,1) $ sufficiently small. The case when $dist(\mathcal{C}_{q_{1}},\mathcal{C}_{q_{2}}) \geq \delta $ easily follows from the continuous inclusions $L^{\infty}(\Omega) \hookrightarrow L^{2}(\Omega) \hookrightarrow H^{-1}(\Omega) $ using the bound $N$ on the norm of the potentials. Therefore we shall henceforth focus on the case $dist(\mathcal{C}_{q_{1}},\mathcal{C}_{q_{2}})< \delta $ where the choice of $\delta $ shall be made clear in the course of the proof.\\
Let us denote $E= \vert \log \ \dist(\mathcal{C}_{q_{1}},\mathcal{C}_{q_{2}}) \vert $ and for $\rho > 0 $ to be chosen later, we set $Z_{\rho}:= \{\xi \in \mathbb{R}^{3}: \vert \xi' \vert < \rho, \vert \xi_{3} \vert < \rho \} $. The integral over the higher frequency modes can be estimated using the bounds on the $L^{2}$-norms of the potentials $q_{1},q_{2} $ and  we can write  
\begin{equation}
\begin{aligned}
\Vert q_{0} \Vert^{2}_{H^{-1}} &=\int_{Z_{\rho}} \frac{\vert \mathcal{F} q_{0}(\xi) \vert^{2}}{1+ \vert \xi \vert^{2}} \ d\xi + \int_{Z^{c}_{\rho}} \frac{\vert \mathcal{F} q_{0}(\xi) \vert^{2}}{1+\vert \xi \vert^{2}}\ d\xi  \leq \int_{Z_{\rho}} \frac{\vert \mathcal{F} q_{0}(\xi) \vert^{2}}{1+ \vert \xi \vert^{2}} \ d\xi  + \frac{C}{\rho^{2}}.
\end{aligned}
\label{estimate-27}
\end{equation}
In order to estimate the integral over the lower frequency modes, that is, the first term in the right hand side of \eqref{estimate-27} we proceed as follows. Using \eqref{estimate-26}, provided $\vert \xi \vert^{2} (\frac{1}{4}+\tau^{2}) > M^{2}+k^{2} $, we can write 
\begin{equation}
\begin{aligned}
\int_{Z_{\rho}} \frac{\vert \mathcal{F} q_{0}(\xi) \vert^{2}}{1+ \vert \xi
\vert^{2}} \ d\xi &\leq C \Big[\int_{Z_{\rho}} \frac{k^{8} e^{4R [\vert \xi
\vert^{2}(\frac{1}{4}+\tau^{2}) - k^{2}]^{\frac{1}{2}}}
\dist(\mathcal{C}_{q_{1}},\mathcal{C}_{q_{2}})^{2} + \epsilon^{2 \alpha} +
\frac{1}{\vert \xi \vert^{2}(\frac{1}{4}+\tau^{2})-k^{2}}}{1+\vert \xi
\vert^{2}} \Big] d\xi+ C \int_{-\rho}^{\rho} \int_{B'(0,\rho)}
\frac{e^{-\frac{\epsilon^{2} (1+ 4 \tau^{2}) \vert \xi' \vert^{2}}{2 \pi}}
}{1+\vert \xi \vert^{2}}d\xi .
\end{aligned}
\label{estimate-28}
\end{equation}
Now we choose $\frac{1}{4}+\tau^{2}=\frac{2 k^{2}+ (\frac{E}{5 R})^{2}}{\vert
\xi \vert^{2}} $. This would imply that $\vert \xi \vert^{2} (\frac{1}{4} +
\tau^{2})-k^{2} = k^{2}+ (\frac{E}{5R})^{2} $ and therefore we shall also have
to choose $E$ such that $(\frac{E}{5R})^{2} > M^{2} $. This, in turn, is
linked to the choice of the $\delta \in (0,1) $. In fact, choosing $\delta $
sufficiently small such that $\dist(\mathcal{C}_{q_{1}},\mathcal{C}_{q_{2}})<
\delta  $, $E$ can be made large enough to fulfil the condition. It will be
worth noting that the choice of $\delta $ depends on the constants $R $ and $M
$ only.  \\
Then 
\begin{equation}
\begin{aligned}
&C \int_{Z_{\rho}} \frac{1}{\vert \xi \vert^{2} (\frac{1}{4}+\tau^{2})-k^{2}}
d\xi \leq \frac{C \rho^{3}}{k^{2}+(\frac{E}{5R})^{2}} \leq \frac{C \rho^{3}}{(k+\frac{E}{5R})^{2}},  \ \text{and} \ 
C \int_{Z_{\rho}} \epsilon^{2 \alpha} d\xi= C \rho^{3} \epsilon^{2 \alpha}.
\end{aligned}
\label{estimate-29}
\end{equation}
Also
\begin{equation}
e^{4R [\vert \xi \vert^{2} (\frac{1}{4}+\tau^{2})-k^{2}]^{\frac{1}{2}}}
\dist(\mathcal{C}_{q_{1}},\mathcal{C}_{q_{2}})^{2}= e^{4R
[k^{2}+(\frac{E}{5R})^{2}]^{\frac{1}{2}}} \dist(\mathcal{C}_{q_{1}},\mathcal{C}_{q_{2}})^{2},
\notag
\end{equation}
and therefore
\begin{equation}
\begin{aligned}
C \int_{Z_{\rho}} k^{8} e^{4R [\vert \xi
\vert^{2}(\frac{1}{4}+\tau^{2})-k^{2}]^{\frac{1}{2}}}
\dist(\mathcal{C}_{q_{1}},\mathcal{C}_{q_{2}})^{2} d\xi &\leq C \rho^{3} k^{8}
e^{4R[k+\frac{E}{5R}]} \dist(\mathcal{C}_{q_{1}},\mathcal{C}_{q_{2}})^{2}.
\end{aligned}
\label{estimate-30}
\end{equation}
To estimate the last term in the right-hand side of \eqref{estimate-28}, we proceed as follows (see also \cite{Choudhury-Heck}).
We note that $1+4 \tau^{2}=\frac{8k^{2}+4(\frac{E}{5R})^{2}}{\vert \xi \vert^{2}} \geq \frac{2[k^{2}+(\frac{E}{5R})^{2}]}{\vert \xi \vert^{2}} \geq \frac{[k+ \frac{E}{5R}]^{2}}{\vert \xi \vert^{2}} $ which implies
\begin{equation}
\begin{aligned}
e^{-\frac{\epsilon^{2} (1+4 \tau^{2}) \vert \xi' \vert^{2}}{2 \pi}} &\leq e^{-\frac{\epsilon^{2} [k+\frac{E}{5R}]^{2} \vert \xi' \vert^{2}}{2 \pi \vert \xi \vert^{2}}}.
\end{aligned}
\notag
\end{equation}
Also since $\vert \xi \vert^{2} < 2 \rho^{2} $, we have $e^{\frac{\epsilon^{2}[k+ \frac{E}{5R}]^{2} \vert \xi' \vert^{2}}{2 \pi \cdot 2 \rho^{2}}} \leq e^{\frac{\epsilon^{2} [k+\frac{E}{5R}]^{2} \vert \xi' \vert^{2}}{2 \pi \vert \xi \vert^{2}}} $ and therefore 
\begin{equation}
e^{-\frac{\epsilon^{2} (1+4 \tau^{2}) \vert \xi' \vert^{2}}{2 \pi}} \leq e^{-\frac{\epsilon^{2} [k+\frac{E}{5R}]^{2} \vert \xi' \vert^{2}}{2 \pi \cdot 2 \rho^{2}}}.
\notag
\end{equation}
Let us choose $\epsilon > 0 $ such that $\epsilon^{2}=\frac{1}{k+\frac{E}{5R}}
$. If required, we can choose $\delta $ smaller again 
such that $\epsilon < \epsilon_{0} $ also holds. Then we can write
\begin{equation}
\begin{aligned}
\int_{-\rho}^{\rho} \int_{B'(0,\rho)} \frac{e^{-\frac{\epsilon^{2}(1+4
\tau^{2}) \vert \xi' \vert^{2}}{2 \pi}}}{1+\vert \xi \vert^{2}} d\xi 
	&\leq \int_{-\rho}^{\rho} \int_{B'(0,\rho)} \frac{e^{-\frac{\epsilon^{2} [k+\frac{E}{5R}]^{2} \vert \xi' \vert^{2}}{4 \pi \rho^{2}}}}{1+\vert \xi \vert^{2}} \ d\xi' d\xi_{3} 
\leq C \rho \int_{0}^{\rho} r e^{-\frac{k+ \frac{E}{5R}}{4 \pi \rho^{2}} r^{2}} \ d\xi' d\xi_{3} \\
&\leq C \rho^{2} \Big[k+\frac{E}{5R} \Big]^{-\frac{1}{2}} \rho \Big[k+\frac{E}{5R} \Big]^{-\frac{1}{2}} \int_{0}^{\infty} u e^{-\frac{1}{4 \pi} u^{2}} \ du \leq C \rho^{3} \Big[k+\frac{E}{5R} \Big].
\end{aligned}
\label{estimate-31}
\end{equation}
We shall now specify our choice of $\rho$. 
Since $\alpha \in (0,1) $, we have
$C \rho^{3} \epsilon^{2 \alpha} = C \rho^{3} \Big(\frac{1}{k+\frac{E}{5R}}\Big)^{\alpha} \geq C \rho^{3} \Big(\frac{1}{k+\frac{E}{5R}}\Big).$\\
Now we choose $\rho > 0 $ such that $\rho^{3}=\Big(k+\frac{E}{5R} \Big)^{\beta} $, where $\beta < \alpha (< 1) $. Then $C \rho^{3} \Big(\frac{1}{k+\frac{E}{5R}} \Big)^{\alpha}= C \Big(\frac{1}{k+\frac{E}{5R}} \Big)^{\alpha-\beta} $ and also $\frac{C}{\rho^{2}}= C. \Big(\frac{1}{k+\frac{E}{5R}} \Big)^{\frac{2 \beta}{3}} $ and $\frac{C \rho^{3}}{(k+\frac{E}{5R})^{2}}= C \Big(\frac{1}{k+\frac{E}{5R}} \Big)^{2-\beta} $.\\
Finally, let $\tilde{\alpha}=\text{min} \{\alpha-\beta, 2-\beta, \frac{2
\beta}{3} \} $. Then from \eqref{estimate-27}--\eqref{estimate-31} we derive
\begin{equation}
\begin{aligned}
\Vert q_{0} \Vert^{2}_{H^{-1}} &\leq C \rho^{3} k^{8} e^{4R \Big[k+\frac{E}{5R} \Big]} \ dist(\mathcal{C}_{q_{1}},\mathcal{C}_{q_{2}})^{2}+ C \Big(\frac{1}{k+\frac{E}{5R}} \Big)^{\tilde{\alpha}} \\
&=C k^{8} \Big(k+\frac{E}{5R} \Big)^{\beta} e^{4R \Big[k+\frac{E}{5R} \Big]} \ dist(\mathcal{C}_{q_{1}},\mathcal{C}_{q_{2}})^{2}+ C \Big(\frac{1}{k+\frac{E}{5R}} \Big)^{\tilde{\alpha}} \\
&\leq C \Big[k^{8} e^{R \Big[k+\frac{E}{5R} \Big]} e^{4R \Big[k+\frac{E}{5R} \Big]} \ dist(\mathcal{C}_{q_{1}},\mathcal{C}_{q_{2}})^{2}+ C \Big(\frac{1}{k+\frac{E}{5R}} \Big)^{\tilde{\alpha}} \Big] \\
&\leq C \Big[k^{8} e^{5Rk} e^{E} \ dist(\mathcal{C}_{q_{1}},\mathcal{C}_{q_{2}})^{2} + \frac{1}{(k+\frac{E}{5R})^{\tilde{\alpha}}}\Big] \\
&\leq C \Big[e^{6Rk} \ dist(\mathcal{C}_{q_{1}},\mathcal{C}_{q_{2}})+ \frac{1}{(k+\frac{E}{5R})^{\tilde{\alpha}}}  \Big].
\end{aligned}
\label{main-estimate-1}
\end{equation}
As already discussed before, the estimate \eqref{main-estimate-1} also holds true when $dist(\mathcal{C}_{q_{1}},\mathcal{C}_{q_{2}}) \geq \delta$. \\
Using \eqref{main-estimate-1}, we can now estimate the $L^{\infty}$-norm of $q_{0}=q_{2}-q_{1}$ by using interpolation. To see this, we recall that given $t_{0}, t, t_{1} $ such that $t_{0}< t_{1} $ and $t= (1-p)t_{0}+p t_{1} $, where $p \in (0,1) $, the $H^{t}$-norm of a function $q$ can be estimated (by the interpolation theorem) as 
\begin{equation}
\Vert q \Vert_{H^{t}} \leq \Vert q \Vert^{1-p}_{H^{t_{0}}} \cdot \Vert  q \Vert^{p}_{H^{t_{1}}}.
\notag
\end{equation}
In our case, we define $\eta > 0 $ such that $s=\frac{3}{2}+ 2 \eta $ and choose $t_{0}=-1, t_{1}=s $ and $t_{2}=\frac{3}{2}+ \eta= s-\eta $. Thus we can write
\begin{equation}
t= (1-p) t_{0}+ p t_{1}, \ \text{where} \ p=\frac{1+s-\eta}{1+s}.
\notag
\end{equation}
Using the Sobolev embedding theorem and the interpolation theorem, we have 
\begin{equation}
\begin{aligned}
\Vert q_{1}-q_{2} \Vert_{L^{\infty}(\Omega)} \leq C \Vert q_{1}-q_{2} \Vert_{H^{\frac{3}{2}+\eta}(\Omega)} &\leq C \Vert q_{1}-q_{2} \Vert^{1-p}_{H^{-1}(\Omega)} \Vert q_{1}-q_{2}\Vert^{p}_{H^{s}(\Omega)} \leq C \Vert q_{1}-q_{2} \Vert^{\frac{\eta}{1+s}}_{H^{-1}(\Omega)} \\
&\leq C \Big(e^{6Rk} \ dist(\mathcal{C}_{q_{1}},\mathcal{C}_{q_{2}})+ \frac{1}{(k+\frac{E}{5R})^{\tilde{\alpha}}} \Big)^{\frac{\eta}{2(1+s)}},
\end{aligned}
\label{main-estimate-2}
\end{equation}
which is the required stability estimate \eqref{main-10}.

\end{document}